\documentclass[10pt]{article}
\pdfoutput=1
\author{Nikola Mirkov, Bosko Rasuo}
\date{}
\usepackage{graphicx}
\usepackage{amsmath}
\usepackage{amsfonts}
\usepackage{listings}
\usepackage{authblk}
\usepackage{placeins}

\title{A Bernstein Polynomial Collocation Method for the Solution of Elliptic Boundary Value Problems}
\author[1]{N. Mirkov\thanks{Email: nmirkov@vinca.rs}}
\author[2]{B. Rasuo}
\affil[1]{University of Belgrade, Institute of Nuclear Sciences 'Vinca', Serbia}
\affil[2]{University of Belgrade, Faculty of Mechanical Engineering, Serbia}

\begin{document}



\maketitle
\begin{abstract}
In this article, a formulation of a point-collocation method in which the unknown function is approximated using global expansion in tensor product Bernstein polynomial basis is presented. Bernstein polynomials used in this study are defined over general interval $[a,b]$. Method incorporates several ideas that enable higher numerical efficiency compared to Bernstein polynomial methods that have been previously presented. The approach is illustrated by a solution of Poisson, Helmholtz and Biharmonic equations with Dirichlet and Neumann type boundary conditions. Comparisons with analytical solutions are given to demonstrate the accuracy  and convergence properties of the current procedure. The method is implemented in an open-source code, and a library for manipulation of Bernstein polynomials \textsf{bernstein-poly}, developed by the authors.
\end{abstract}

\section{Introduction}
Bernstein polynomials are known because of their many useful properties \cite{Lorentz}. When restricted to the unit interval they are used in the definition of B\'{e}zier curves, which are very important tools in computer graphics. On their own, they have enjoyed increased attention in recent years, specifically as means to represent solutions to differential equations. Bhatti and Bracken \cite{Bhatti2007} present a Galerkin method which uses Bernstein polynomials as trial functions for solution of two-point boundary value problem, Yousefi et. al \cite{Yousefi} use Bernstein polynomials in Ritz-Galerkin method to approximate solution of hyperbolic PDE with an integral condition. Doha et. al. \cite{Doha} present the solution method for high-even-order ordinary differential equations. \\
Unlike all listed examples, which use Bernstein polynomials in some sort of Galerkin method, in this article we develope a collocation method. The motivation to do so comes from the nature of Bernstein polynomials themselves. When function is expanded in Bernstein polynomial basis, the expansion coefficients also respresent the nodal values of the expanded function. For two-point boundary value problems this leads to very simple implementation, allowing direct imposition of the boundary conditions at the end-points of the interval. 
When the unknown function is approximated as a global expansion in Bernstein polynomials basis, and used in a point-collocation method, the approach becomes very similar to pseudospectral methods \cite{Boyd}, \cite{Fornberg} \cite{ICASEref}. Unlike their polynomial counterparts used in pseudospectral methods (e.g. Chebyshev and Legandre), Bernstein polynomials are not orthogonal, and many examples of basis transformation exist when this is necessary (eg. \cite{transformation}). \\
For the present article we assume that $L$, the elliptic differential operator in interior, and $B$, the boundary operator are linear, and that they define well posed boundary value problem
\begin{equation} \label{eq:L-operator}
L[u] (\mathbf{x}) =  f(\mathbf{x}), \;\;\; \mathbf{x} \in \Omega.
\end{equation}
\begin{equation} \label{eq:B-operator}
B[u] (\mathbf{x}) =  g(\mathbf{x}), \;\;\; \mathbf{x} \in \partial \Omega.
\end{equation}
Where $\Omega \in \mathbb R^2 $ is a given rectangular domain.\\
In the remainder of the paper we give a brief review of Bernstein polynomial properties, relevant to the proposed method. The solutions to PDEs defined over two-dimensional domains can be represented by surfaces embedded in three-dimensional Eucledian space. This serves as rationale of Section 3, in which we describe surfaces defined by expansions in tensor product Bernstein polynomial basis, also giving expressions for elliptic operators that we will use subsequently. Then, in Section 4, we describe the fomulation of the proposed collocation method. Finally we give couple of example solutions, as well as their error, and convergence analysis in Section 5.

\section{Properties of Bernstein Polynomials} \label{sec:b-poly}
In further discussion we consider only generalized Bernstein polynomials, those defined over arbitrary interval $[a,b]$, and simply call them Bernstein polynomials.\\ 
The Bernstein polynomials of $n$\textit{th} degree form a complete basis over $[a,b]$, and are defined by
\begin{equation}
B_{i,n}(x)= \binom{n}{i}\frac{(x-a)^i(b-x)^n-i}{(b-a)^n}, \;\;\; i=0,1,...,n,
\end{equation}
where binomial coefficients are given by 
\begin{equation} \label{eq:binom-frac}
\binom{n}{i} = \frac{n!}{i!(n-i)!}.
\end{equation}
 
They satisfy symmetry $B_{i,n}(x) = B_{n-i,n}(1-x)$, positivity $B_{i,n}(x) \geq 0$, and they form partition of unity $ \sum_{i=0}^n B_{i,n}(x) = 1$ on defining interval $x \in [a,b]$. For $i=0$, and $i=n$, they have value equalt to one at $x=a$, and $x=b$, respectively. Otherwise, they have a unique local maximum occuring at $x=i/n$, having the value $B_{i,n}(i/n) = i^i n^{-n} (n-i)^{n-i} \binom{n}{i}$.\\  
In a recent article \cite{Doha}, authors present explicit formula for calculation of arbitrary order derivatives of Bernstein basis functions of any order defined over standard interval $[0,1]$. In an atempt to formulate collocation method for solution of BVP's over arbitary intervals, first step would be to write down this expression for Bernstein polynomials defined for arbitrary interval.

Let $B_{i,n}(x)$ define $i-th$ basis function of $n-th$ order, at point $x$, where $x \in [a,b]$, and $i \in [0,n]$. Then, the derivative of order $p$ of such basis function can be expressed as:
\[ 
D^pB_{i,n}(x) = \frac{n!}{(n-p!)(b-a)^p} \sum_{k=max(0,i+p-n)}^{min(i,p)} \binom p k  B_{i-k,n-p}(x) 
\]

Useful proprerty of Bernstein basis functions is that they all vanish at endpoints of the interval, except the first and the last one, which are equal to one at $x=a$ and $x=b$ respectively. 

The factorial formula Eq. \ref{eq:binom-frac} for binomial coefficients is known not to be most numerically efficient for large numbers, and the alternative representations exist, most numerically efficient being the multiplicative formula, which we will use here
\begin{equation}
\binom{n}{i}  = \prod_{i=1}^k \frac{n-(k-i)}{i}.
\end{equation}

\section{Surfaces defined by tensor product of Bernstein polynomials} 
\label{sec:tensor-product}
Following equations will help us in developing a collocation method using Bernstein polynomials for boundary value problems defined over two-dimensional domains.

Let the following equation represent a surface in $\mathbb R^3$ defined by a tensor product of Bernstein basis functions

\begin{equation} \label{eq:interpolant}
f(x,y) = \sum_{i=0}^n \sum_{j=0}^m \beta_{i,j} B_{i,n}(x) B_{j,m}(y).
\end{equation}

We can define p-th order partial derivative in $x$ -axis direction in a following way

\begin{equation} \label{eq:derx}
\frac{\partial^p f(x,y)}{\partial x^p} = \sum_{i=0}^n \sum_{j=0}^m \beta_{i,j} D^p B_{i,n}(x) B_{j,m}(y),
\end{equation}

and an q-th order derivative in $y$ -axis direction

\begin{equation} \label{eq:dery}
\frac{\partial ^q f(x,y)}{\partial y^q} = \sum_{i=0}^n \sum_{j=0}^m \beta_{i,j} B_{i,n}(x) D^q B_{j,m}(y).
\end{equation}

Finally mixed derivative of order $p+q$ is defined as follows:

\begin{equation}
\frac{\partial ^{p+q} f(x,y)}{\partial x^p y^q} = \sum_{i=0}^n \sum_{j=0}^m \beta_{i,j} D^p B_{i,n}(x) D^q B_{j,m}(y).
\end{equation}
This paper deals with elliptic equations, therefore we will define two differential operators, most common in second and forth order boundary value problems.\\
Laplacian operator can be written as

\[  \Delta f(x,y) = \frac{\partial^2 f(x,y)}{\partial x^2 } + \frac{\partial^2 f(x,y)}{\partial y^2} 
= \]
 
\begin{equation}
 =  \sum_{i=0}^n \sum_{j=0}^m \beta_{i,j} B_{j,m} (y) D^2 B_{i,n}(x) + \sum_{i=0}^n \sum_{j=0}^m \beta_{i,j} B_{i,n} (x) D^2 B_{j,m},
\end{equation}

which can further be simplified to

\begin{equation}
\Delta f(x,y) =  \sum_{i=0}^n \sum_{j=0}^m \beta_{i,j} \left[ B_{j,m} (y) D^2 B_{i,n}(x) + B_{i,n} (x) D^2 B_{j,m}(y) \right].
\end{equation}

Bicharmonic operator or bilaplacian can be written as 

\[ \Delta^2 f(x,y) = \frac{\partial^4 f(x,y)}{\partial x^4 } +2 \frac{\partial^4 f(x,y)}{\partial x^2 \partial y^2 } + \frac{\partial^4 f(x,y)}{\partial y^4} = \]

\begin{equation}
= \sum_{i=0}^n \sum_{j=0}^m \beta_{i,j} \left[ B_{j,m} (y) D^4 B_{i,n}(x)+  2 D^2 B_{i,n} (x) D^2 B_{j,m}(y)  + B_{i,n} (x) D^4 B_{j,m}(y) \right].
\end{equation}
 
\section{Collocation method formulation}
In present method, it is assumed that a variable can be expressed as an approximation in the form of global expansion in tensor product Bernstein polynomial basis. \\
We look for the approximate solution in the form (\ref{eq:interpolant}) for $(x,y) \in \bar{\Omega} = \Omega \cup \partial \Omega$. The expansion coefficients $\beta_{i,j}$ are unknown, and need to be determined.  \\
The collocation method, that we use to find the unknown coefficients, is a numerical procedure in which we require that solution satisfies differential equation exactly in a discrete set of points, known as collocation points.
Number of collocation points has to match the number of unknowns.\\
If the order of Bernstein polynomial is $n$ and $m$ in $x$ and $y$-axis directions respectively, then there is $n_c = (n+1)\times(m+1)$ unknown coefficients in global polynomial expansion. If we are solving homogenous Dirichlet problem the number of unknown is $n_{cd} = (n-1)\times(m-1)$ . Reason for this is a property of Bernstein polynomials that we take the advantage of, namely the property to vanish at end points, except for the first and the last one, that we mentioned in Section \ref{sec:b-poly}. If the function value at the boundary is equal to zero, we eliminate the non-vanishing basis functions by setting their associated expansion coefficients to zero, and continue by solving the problem for the smaller number of unknowns $n_{cd}$.\\ 
Discretized equation set is obtained by substituting the unknown function  and it's derivatives with the appropriate representation in Bernstein polyniomial basis, as described in Section \label{sec:tensor-product}. 
Substituting these Bernstein polynomial interpolants into the original equation set, and applying it at the set of collocation points, one arrives at a linear system of equations
\begin{equation}\label{eq:linsys}
A b = c.
\end{equation}
An element of the system matrix $a_{i,j}$ is obtained by evaluating the differential operator on the tensor product of $i$\textit{th} and $j$\textit{th} Bernstein polynomial basis function, at the collocation point defined by Cartesian coordinate pair $\left( x_i, y_i \right)$.

\begin{equation} \label{eq:a_assemble}
a_{i,j} = L \left [ B_{i,n}(x_i) B_{j,m}(y_i) \right ], \;\;\; i = 1,..,n-1, j = 1,..,m-1.
\end{equation}

An element of the right-hand side vector $c_i$ is an evaluation of the right-hand side function at the $i$\textit{th} collocation point defined by Cartesian coordinates $\left( x_i, y_i \right)$.

\begin{equation}
c_i = f(x_i,y_i), \;\;\; i=1,...,n_{cd}.
\end{equation}

For the cases other than those with homogenous Dirichlet boundary conditions, the number of unknowns, as being said, is larger and linear system has the additional rows defined by

\begin{equation} \label{eq:b_assemble}
a_{i,j} = B \left [ B_{i,n}(x_i) B_{j,m}(y_i) \right ], 
\end{equation}
\begin{equation}
c_i = g(x_i,y_i), 
\end{equation}
Where coordinates $\left( x_i, y_i \right)$ define colocation points that belong to the domain boundary. Index pairs $(i,j)$ belong to the set $\left(0,1...m-1 \right) \cup \left(n,1...m-1 \right) \cup \left(1...n-1,0 \right) \cup \left( 1...n-1,m \right) \cup \left( 0,0\right) \cup \left( n,0\right) \cup \left( 0,m\right) \cup \left( n,m\right) $. There is in total $2(n+m)$ additional rows required to specify boundary conditions.  \\
The change in boundary conditions (\ref{eq:B-operator}), amounts to changing a few rows in the matrix $A$, as well as in the right-hand side vector $c$. Using expressions defined in Section 3. we may define any type of boundary conditions. In particular, the operator $B$ takes the form of Eq.(\ref{eq:interpolant}) when we need to impose Dirichlet boundary conditions, or (\ref{eq:derx}),(\ref{eq:dery}), with $p=q=1$ for Neumann boundary conditions.\\
Solution vector $b$ is as an auxilliary one-parameter array, it's values are moved to two-parameter array of Bernstein polynomial expansion coefficient $\beta$. Mapping has the following form
\begin{equation} \label{eq:cmap}
\beta_{i,j} = b_{inode},
\end{equation}
Where $inode = (m+1)*(i-1)+j$, and $i =0,...,n, j=0,...,m$. If only inner expansion coefficients are sought for, in the case of homogenous Dirichlet problems, coefficient mapping takes the same form (\ref{eq:cmap}), but the indexes are defined by $inode = (m-1)*(i-1)+j-1$, and $i=1,...,n-1, j=1,...,m-1$.

Finally we note that the distribution of points on a tensor product grid may be uniform, or non-uniform depending on situation. The full linear system \ref{eq:linsys} is solved using LU decomposition with partial pivoting. 

\section{Example Application}
Collocation method formulated in previous section is implemented in \textsf{bernstein-poly}, a library for manipulation of Bernstein polynomials. It is an open-source code developed by the authors, written in Python \cite{b-poly-code}. In what follows, we will present couple of examples, with the purpose to illustrate  validity and accuracy of the present collocation method. In all examples we use the same order of Bernstein polynomials in both $x$ and $y$ axis direction.

\subsection*{Example 1.}  Consider Poisson equation

\begin{equation}
\label{eq:poisson} 
\Delta u(x,y) = f(x,y), \;\;\; x,y \in [-1,1],
\end{equation}
with homogenous Dirichlet boundary conditions. Source term is defined as $f(x,y) = -2 \pi^2 sin(\pi x) \, sin(\pi y) $.
Exact solution for the given problem is
\begin{equation}
u_{exact} = sin(\pi x) \, sin(\pi y)
\end{equation}
Fig. 1 shows numerical solution for approximation using Bernstein polynomials of order n=21. The absolute numerical solution error $abserr(i,j) = \Vert u(i,j) - u_{i,j}^e \Vert $ , where $u_{i,j}^e$ represents the exact solution, is estimated at each $(i,j)$th collocation node, and the result is presented in Fig. 2. To study how accuracy changes with incresing the order of polynomial approximation we use $L^2$ relative error norm, defined in the following way
\begin{equation}
L^2 Error = \frac{\sum_{i,j=0}^{n,m} \left( u_{i,j} - u_{i,j}^e \right)^2 }{ \sum_{i,j=0}^{n,m} (u_{i,j}^e)^2}.
\end{equation}
Results of the analysis of $L^2$ relative error norm are sumarized in Table 1. To study order of accuracy it is useful to plot $L^2$ relative error norm as a function of polynomial order, $n$, which is done in \textit{log-log} plot in Fig. \ref{fig:ex1l2err}. We observe exponential decay in the error for polynomial orders up to $n=17$, after which the error decrease slows down. Eventually for order of Bernstein basis polynomials higher than $n=21$, after which the error continually grows.

\begin{figure}[htb]
\centering
\includegraphics[width=0.8\textwidth]{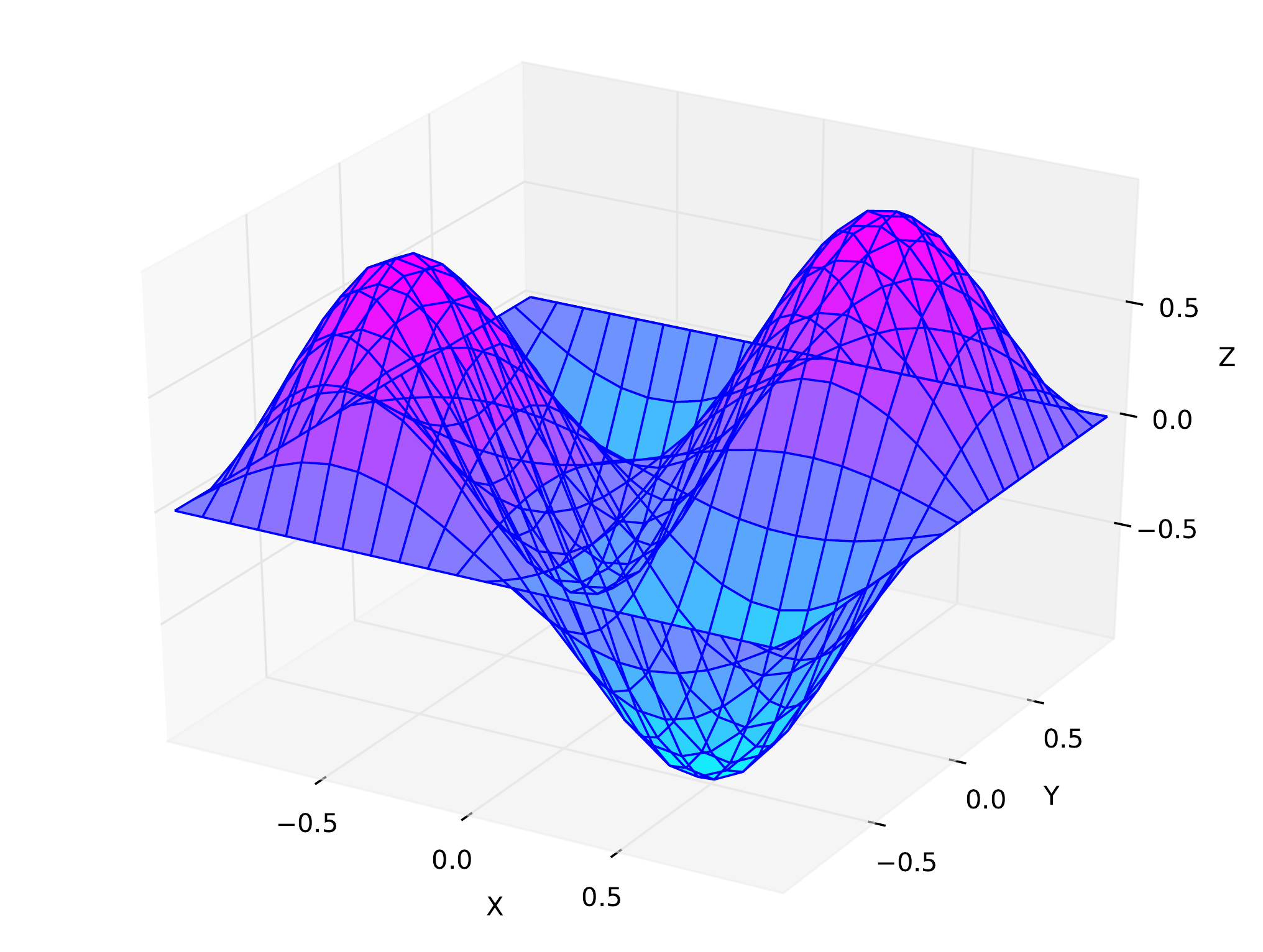}
\caption{Bernstein polynomial collocation method solution of Poisson equation (Example 1), n=20.}
\label{fig:ex1}
\end{figure}

\begin{figure}[htb]
\centering
\includegraphics[width=1.0\textwidth]{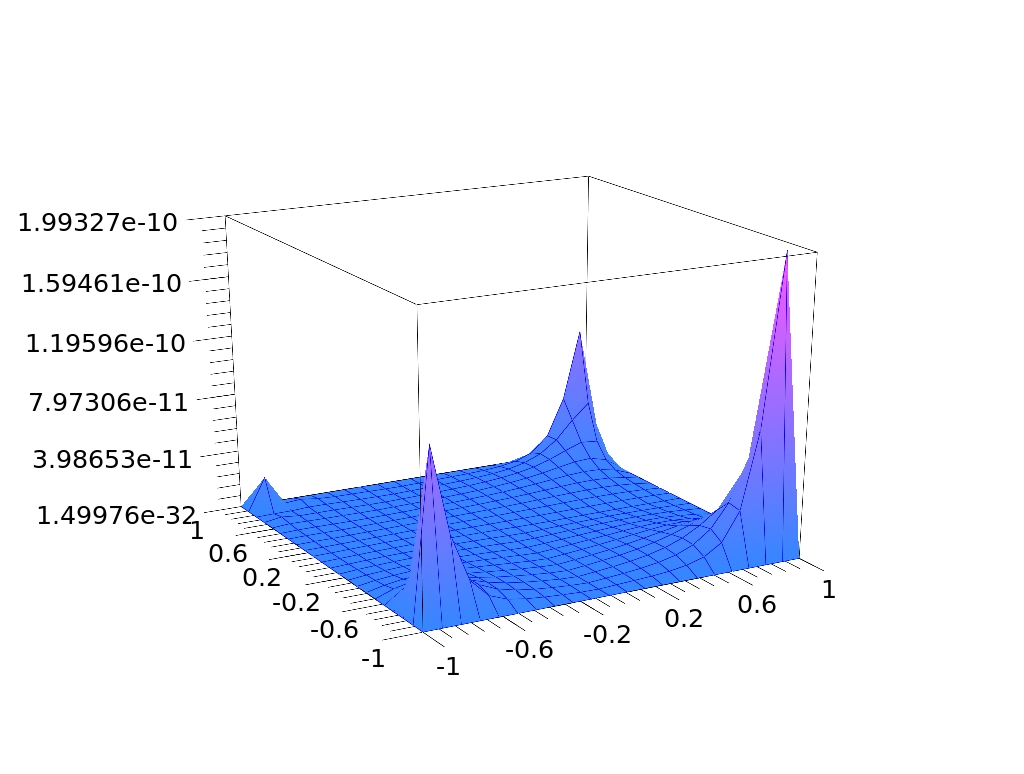}
\caption{Absolute error (Example 1), n=20.}
\label{fig:ex1abserr}
\end{figure}

\begin{figure}[htb]
\centering
\includegraphics[width=0.7\textwidth]{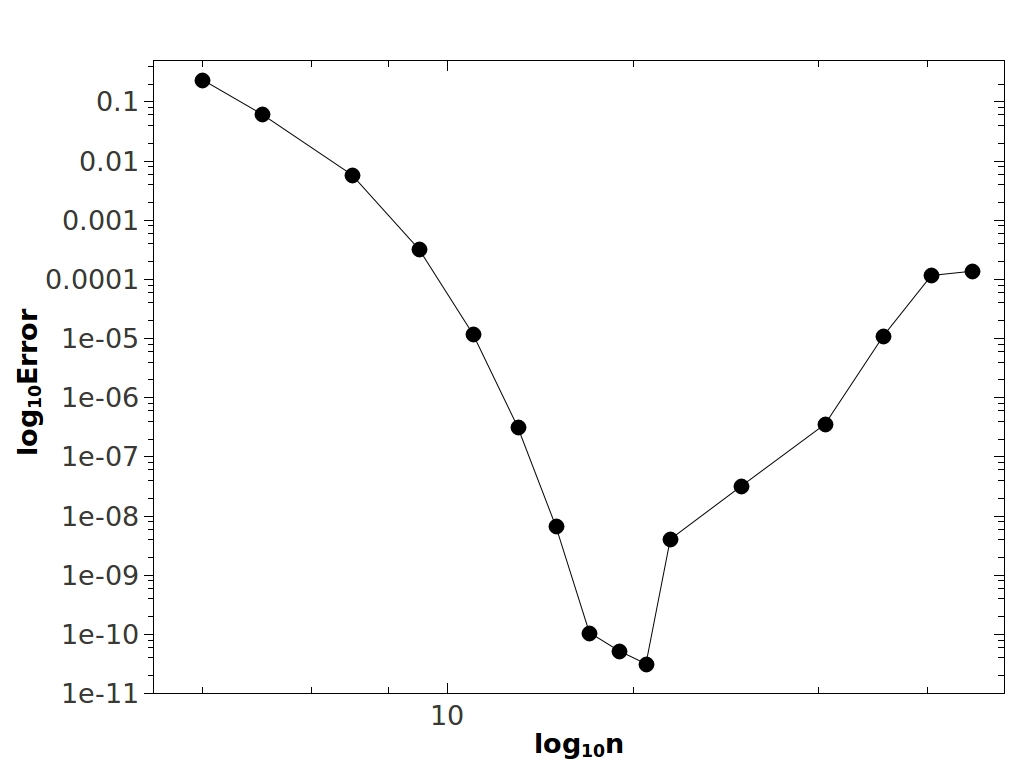}
\caption{Example 1. $L^2$ relative error norm.}
\label{fig:ex1l2err}
\end{figure}


\begin{table}
\centering
\begin{tabular} {|c|c|}
\hline
n & $L^2$ rel. error \\ \hline \hline
11 & $1.171 \times 10^{-5}$  \\
13 & $3.170 \times 10^{-7}$  \\
15 & $6.536 \times 10^{-9}$  \\
17 & $1.049 \times 10^{-10}$ \\
19 & $5.137 \times 10^{-11}$ \\
21 & $3.111 \times 10^{-11}$ \\
23 & $3.966 \times 10^{-9}$  \\
30 & $3.146 \times 10^{-08}$ \\
41 & $3.499 \times 10^{-07}$ \\
51 & $1.091 \times 10^{-05}$ \\
61 & $1.140 \times 10^{-04}$ \\
71 & $1.346 \times 10^{-04}$ \\
\hline
\end{tabular}
\caption{$L^2$ relative error norm for the Example 1.}
\label{l2err1}
\end{table}

\subsection*{Example 2.}  
\begin{figure}[htb]
\centering
\includegraphics[width=0.8\textwidth]{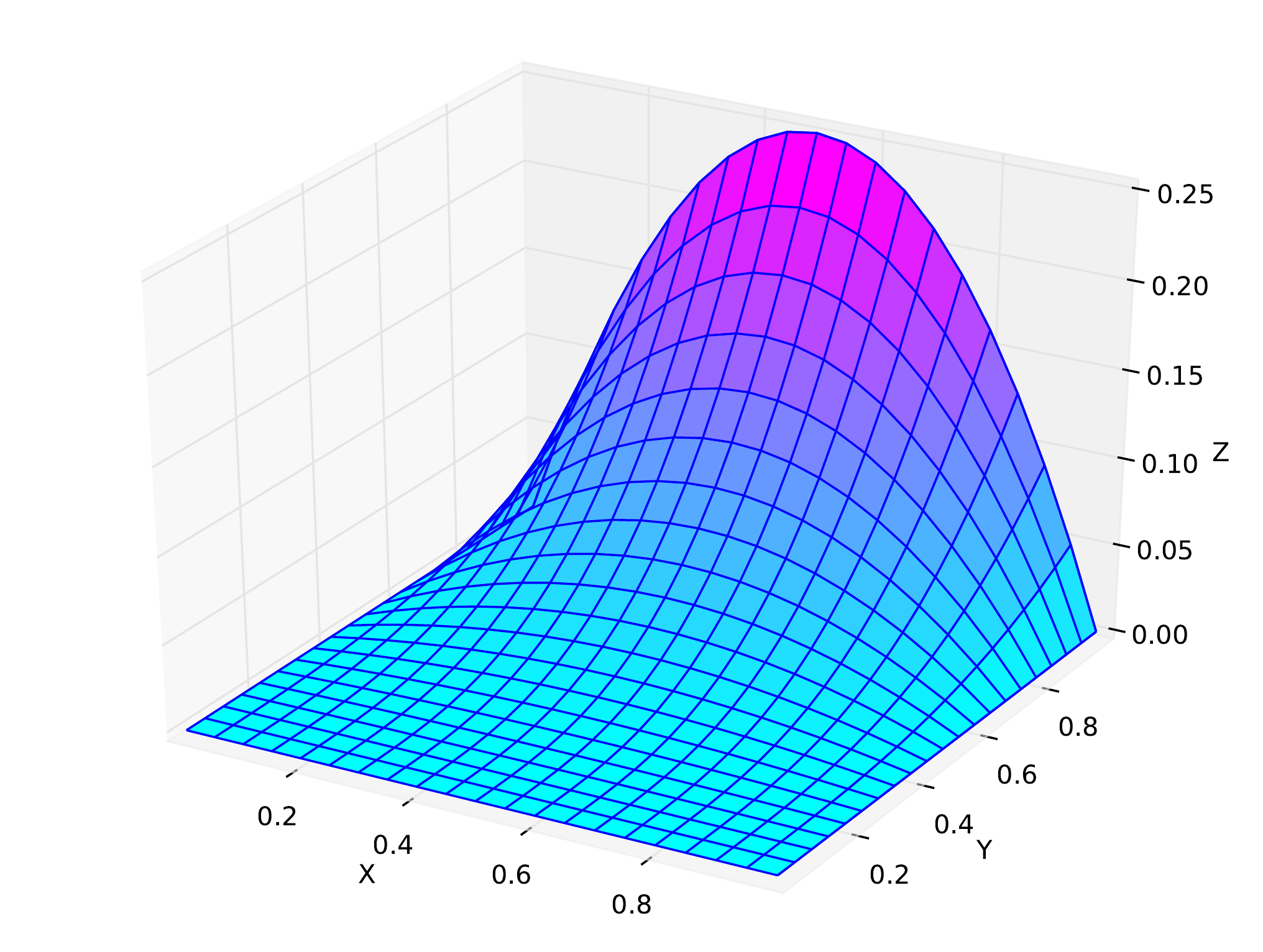}
\caption{Solution of Poisson equation, non-homogenous Dirichlet BCs (Example 2), n=20.}
\label{fig:ex2}
\end{figure}

\begin{figure}[htb]
\centering
\includegraphics[width=0.8\textwidth]{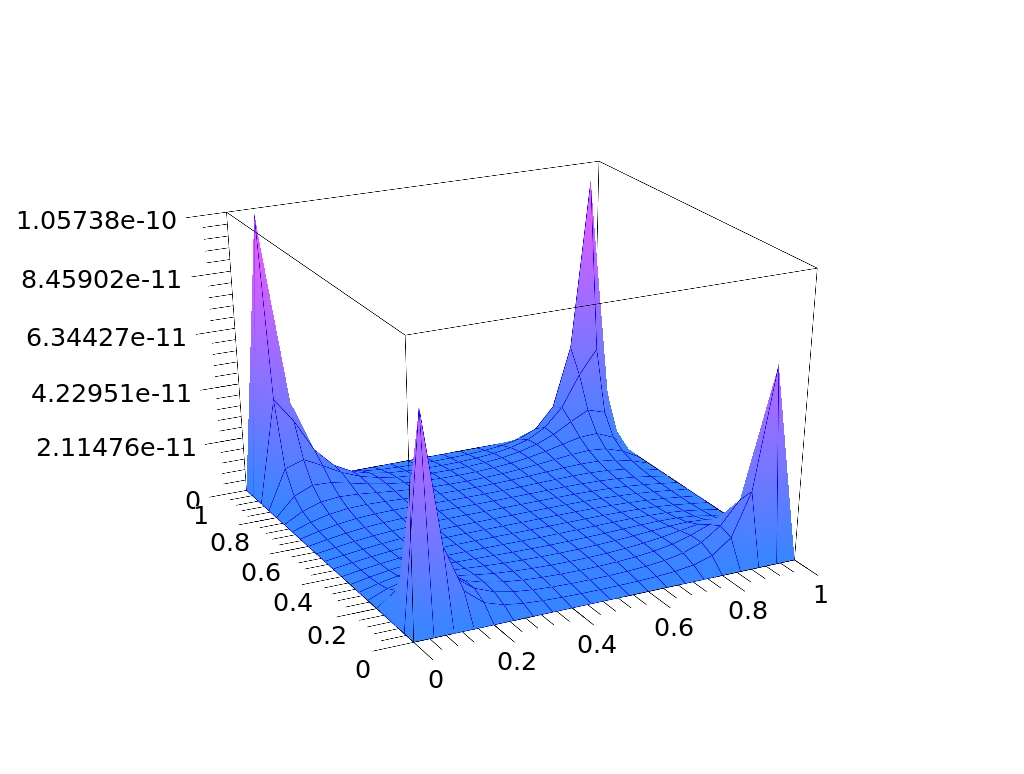}
\caption{Absolute error (Example 2), n=20.}
\label{fig:ex2abserr}
\end{figure}

In this example we consider the Poisson equation with non-homogenous Dirichlet boundary conditions. Treatment of this problem differs from the previous one algorithmically in a way we impose boundary conditions. For every collocation point at the boundary we write an additional linear equation and we solve linear system of size $(n_c)\times(n_c)$.\\
Source term in the Poisson equation (\ref{eq:poisson}) for this problem is defined by $f(x,y) = 6xy(1-y)-2x^3$, and solution domain by unit square $x,y \in [0,1]$. The exact solution for this problem is
\begin{equation}
u_{exact} = y(1-y)x^3
\end{equation}
The figures \ref{fig:ex2}, and \ref{fig:ex2abserr} show numerical solution and absolute error distribution for the case of $n=20$. very high accuracy is achieved quite "early", for $n=12$ it reaches order of $10^{-15}$, as seen in Table \ref{l2err2}.

\begin{table}
\centering
\begin{tabular} {|c|c|}
\hline
n & $L^2$ rel. error \\ \hline \hline
12 & $5.841 \times 10^{-15}$ \\
14 & $2.595 \times 10^{-14}$  \\
16 & $1.754 \times 10^{-13}$ \\
18 & $5.039 \times 10^{-12}$ \\
20 & $1.544 \times 10^{-10}$ \\
30 & $1.907 \times 10^{-8}$ \\
\hline
\end{tabular}
\caption{$L^2$ relative error norm for the Example 2.}
\label{l2err2}
\end{table}

\subsection*{Example 3.} 
\begin{figure}
\centering
\includegraphics[width=0.8\textwidth]{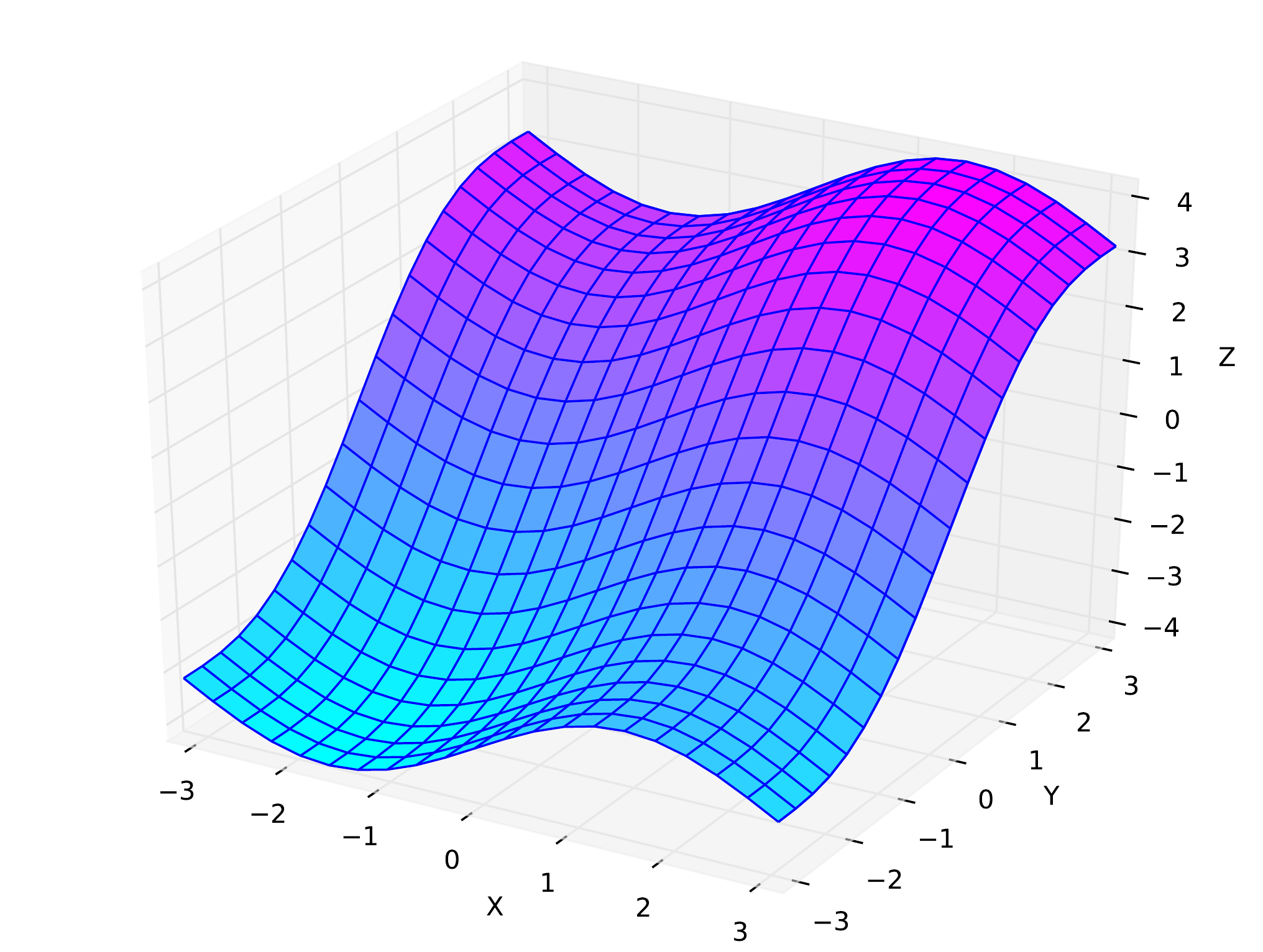}
\caption{Solution of Helmholtz equation, non-homogenous Dirichlet BCs (Example 3), n=20.}
\label{fig:ex3}
\end{figure}

\begin{figure}
\centering
\includegraphics[width=0.8\textwidth]{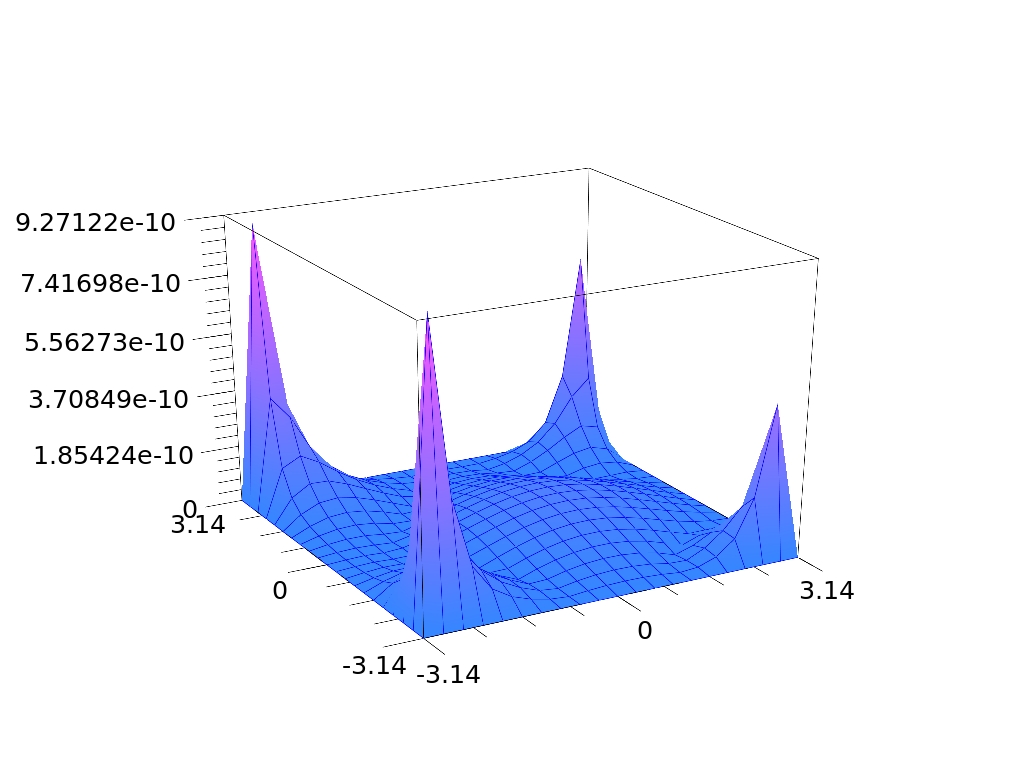}
\caption{Absolute error (Example 3), n=20.}
\label{fig:ex3abserr}
\end{figure}

Next example problem is defined by the Helmholtz equation 
\begin{equation}
(\Delta + \lambda)u = f(x,y).
\end{equation}
The problem is originally found in \cite{Hon}, and is defined by $\lambda = 1$, $f(x,y) = x$, and non-homogenous Dirichlet boundary conditions which are derived from the exact solution. Domain of solution is square $x,y \in [-\pi, \pi]$. Using Bernstein polynomials defined over general interval proves practical for this problem, as mapping of the domain to the unit square is not necessary. \\
This problem is exactly solvable and the solution is 
\begin{equation}
u_{exact} =  sin(x)+sin(y)+x.
\end{equation}
As in previous examples we show numerical solution Fig. \ref{fig:ex3} and the distribution of the absolute error Fig. \ref{fig:ex3abserr}. Table \ref{l2err3} lists $L^2$ relative error norm variation with increasing order of polynomial basis.

\begin{table}
\centering
\begin{tabular} {|c|c|}
\hline
n & $L^2$ rel. error \\ \hline \hline
12 & $9.035 \times 10^{-6}$ \\ 
14 & $2.430 \times 10^{-7}$ \\
16 & $4.992 \times 10^{-9}$ \\
18 & $8.107 \times 10^{-11}$ \\
20 & $4.057 \times 10^{-11}$ \\
22 & $1.051 \times 10^{-9}$ \\
30 & $1.533 \times 10^{-7}$ \\
\hline
\end{tabular}
\caption{$L^2$ relative error norm for the Example 3.}
\label{l2err3}
\end{table}

\subsection*{Example 4.} 
\begin{figure}
\centering
\includegraphics[width=0.8\textwidth]{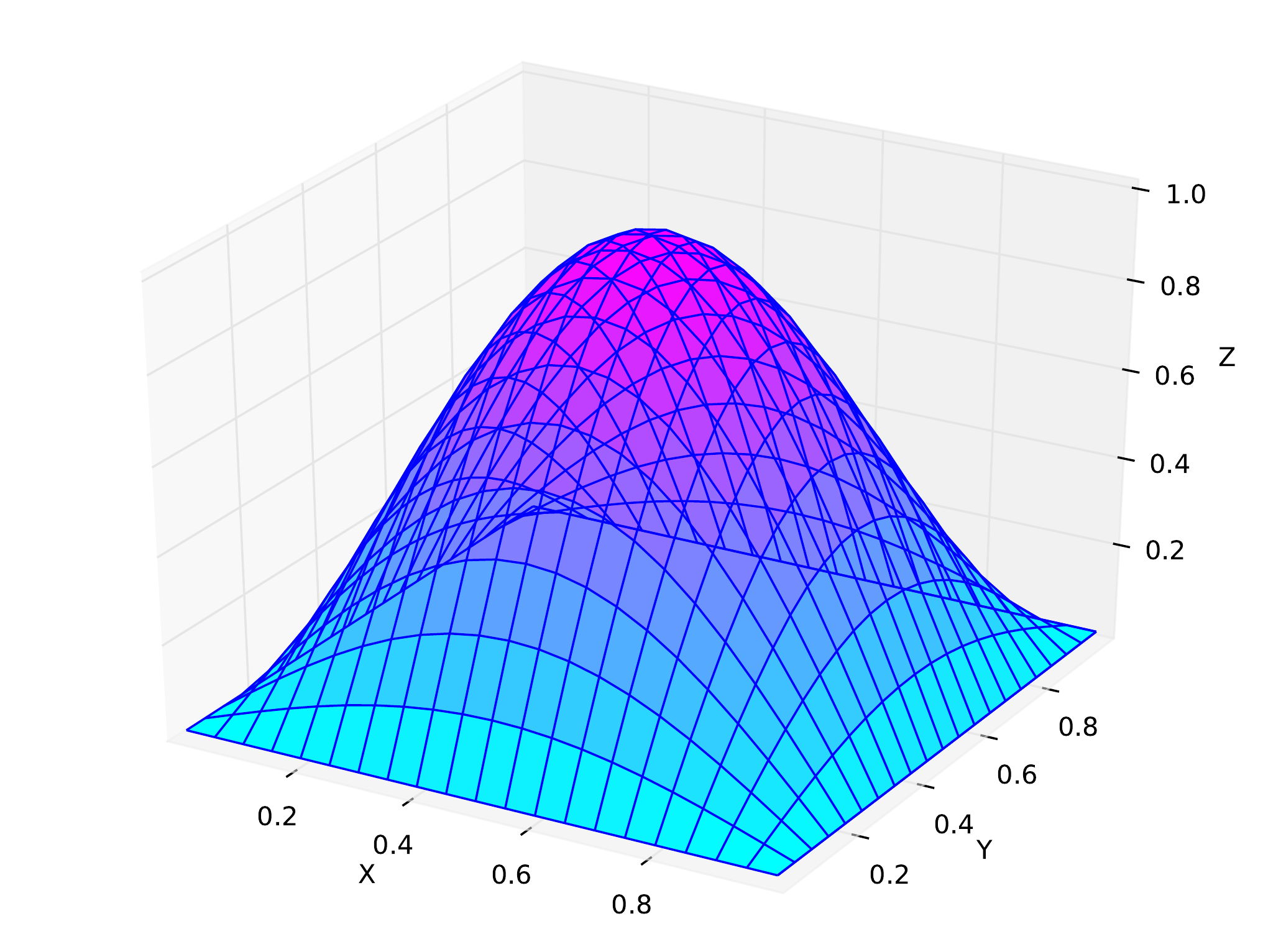}
\caption{Solution of Biharmonic equation (Example 4), Type I BCs, n=20.}
\label{fig:ex4}
\end{figure}

\begin{figure}
\centering
\includegraphics[width=0.8\textwidth]{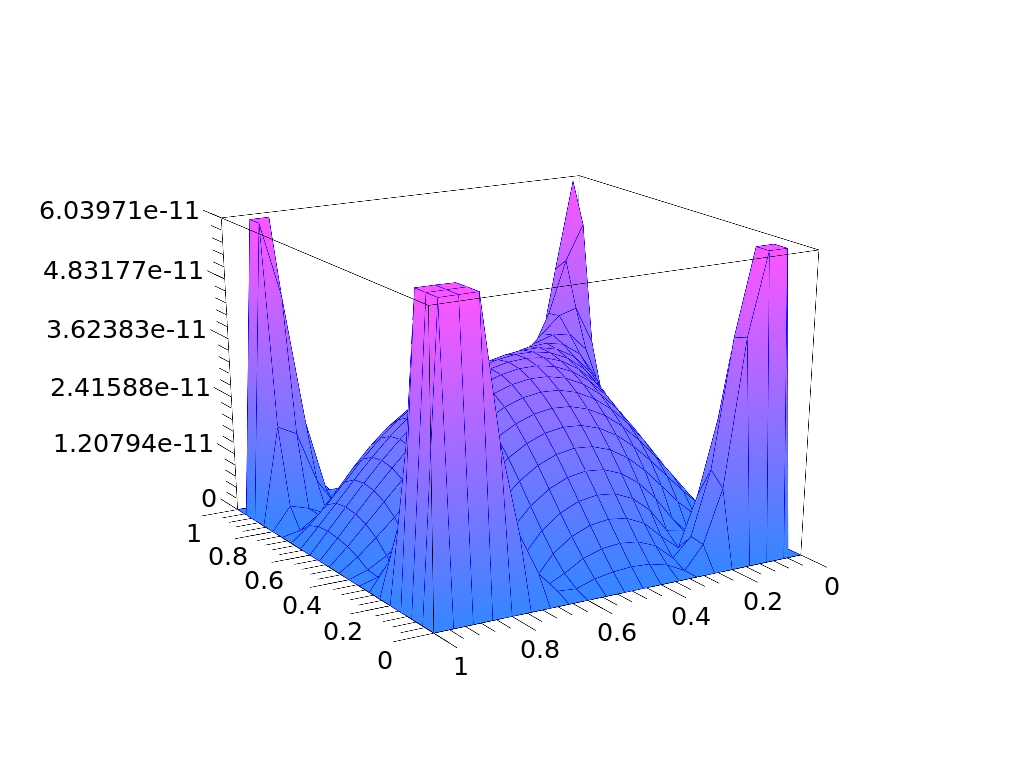}
\caption{Absolute error (Example 4), n=20.}
\label{fig:ex4abserr}
\end{figure}

In next example we give the solution to biharmonic equation, which is often encountered in the theory of ellasticity, as it describes deflections of loaded plates. 

\begin{equation}
\Delta^2 u(x,y) = f(x,y).
\end{equation}
Two types of boundary conditions can be defined for this problem (I) $u$ and $\partial^2 u/ \partial n^2$ or (II) $u$ and $\partial u/ \partial n$. Both cases are interesting on their own because they require different solution approach. We will split discussion on biharmonic equation in two parts, present Example we will treat Type I and in the next one Type II boundary conditions.  \\
In the case of Type I boundary conditions, biharmonic equation can be split into two coupled Poisson equations
\begin{equation}
\Delta v(x,y) = f(x,y),
\end{equation}
\begin{equation}
\Delta u(x,y) = v(x,y).
\end{equation}
The first example solution of biharmonic equation, taken from \cite{Yakhot}, will deal with the case of simply supported rectangular plate ${0 \leq x \leq a, 0 \leq x \leq b}$. Homogenous boundary conditions for both the function value and it's second derivatives in the direction normal to boundary are prescribed. \\
Source function, that describes load distribution over the surface in theory of plates is given by
\begin{equation}
f(x,y) = \pi^4 \left( \frac{m^2}{a^2} +  \frac{n^2}{b^2} \right)^2 sin \frac{m\pi x}{a} sin \frac{n\pi y}{b}.
\end{equation}
Biharmonic equation defined in such a way, allows the exact solution
\begin{equation}
u_{exact}(x,y) =sin \frac{m\pi x}{a}sin \frac{n\pi y}{b}.
\end{equation}
For our purpose we set values $m=n=1$, $a=b=1$. Example solution for the order of polynomial, $n=20$ is shown in Fig \ref{fig:ex4}. The peaks in the absolute error near the corner nodes, as noticed in previous examples, is inherent   to the present method, therefore we set an upper limit on the vertical axis here, to be able to get better picture of the distribution of absolute error in the rest of domain, Fig. \ref{fig:ex4abserr}. Table \ref{l2err4} shows variation in the error norm, in which the trend conforms to the one in previous examples. 

\begin{table}
\centering
\begin{tabular} {|c|c|}
\hline
n & $L^2$ rel. error \\ \hline \hline
10 & $6.538 \times 10^{-8}$ \\
12 & $4.854 \times 10^{-10}$ \\
14 & $2.832 \times 10^{-12}$ \\
16 & $1.284 \times 10^{-12}$\\
20 & $9.467 \times 10^{-11}$ \\
30 & $5.939 \times 10^{-8}$\\
\hline
\end{tabular}
\caption{$L^2$ relative error norm for the Example 4.}
\label{l2err4}
\end{table}

\subsection*{Example 5} 
\begin{figure}
\centering
\includegraphics[width=0.8\textwidth]{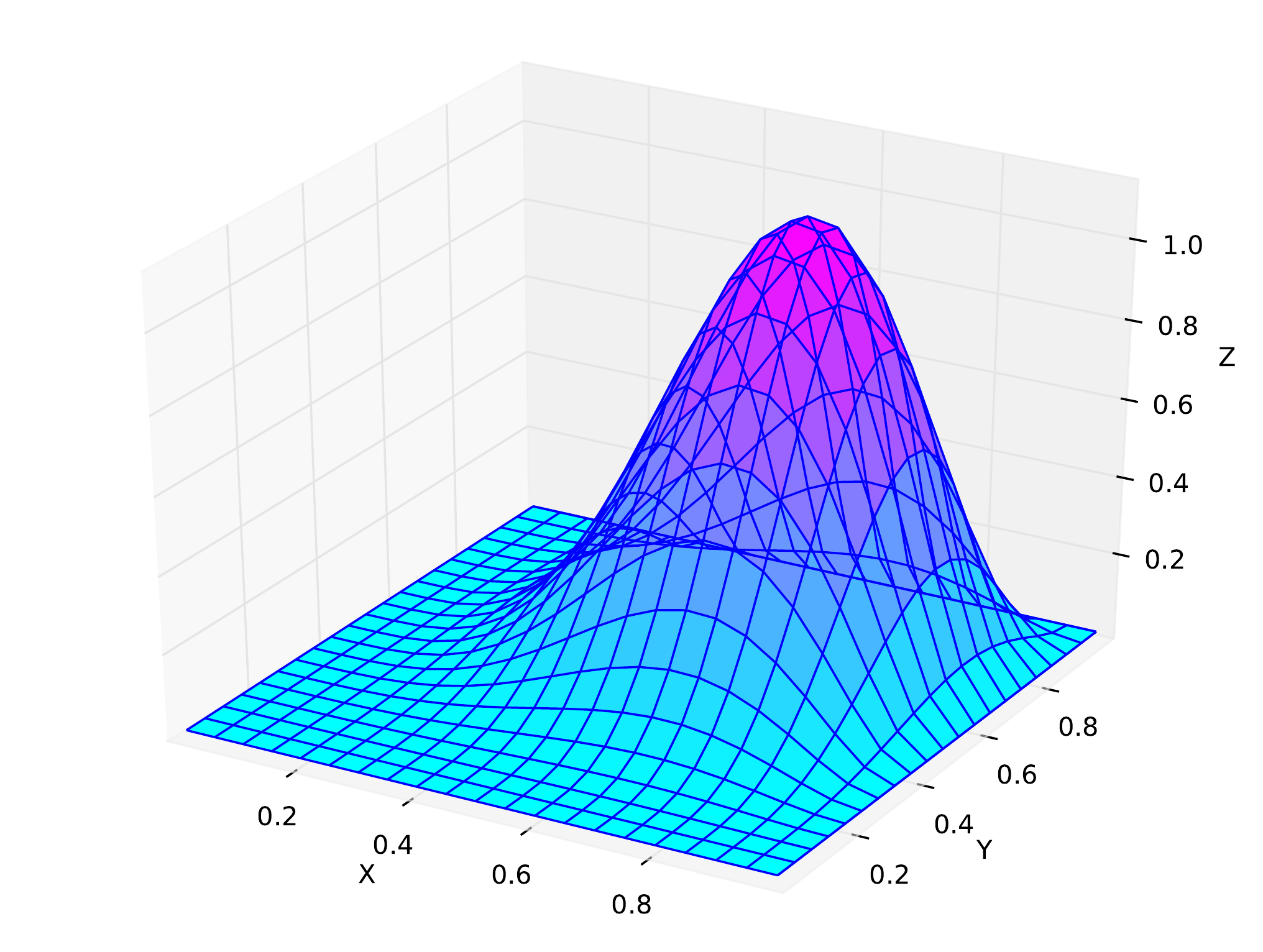}
\caption{Solution of Biharmonic equation (Example 5a), Type II BCs, n=20.}
\label{fig:ex5a}
\end{figure}

\begin{figure}
\centering
\includegraphics[width=0.8\textwidth]{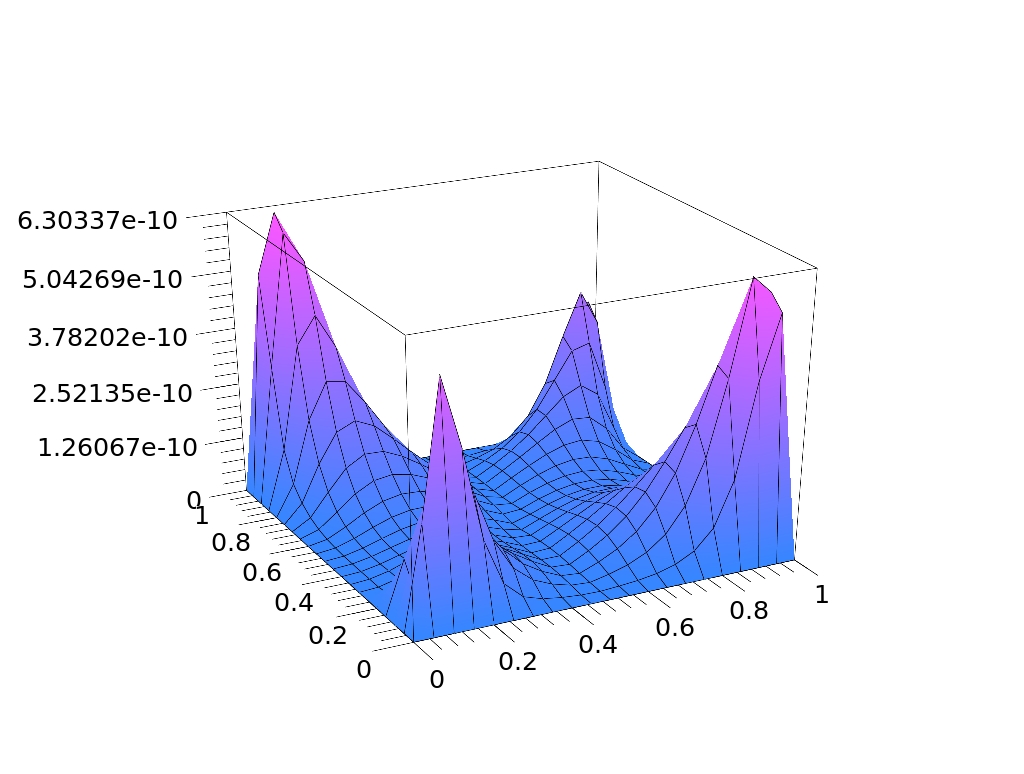}
\caption{Absolute error (Example 5a), n=20.}
\label{fig:ex5aabserr}
\end{figure}

\begin{figure}
\centering
\includegraphics[width=0.8\textwidth]{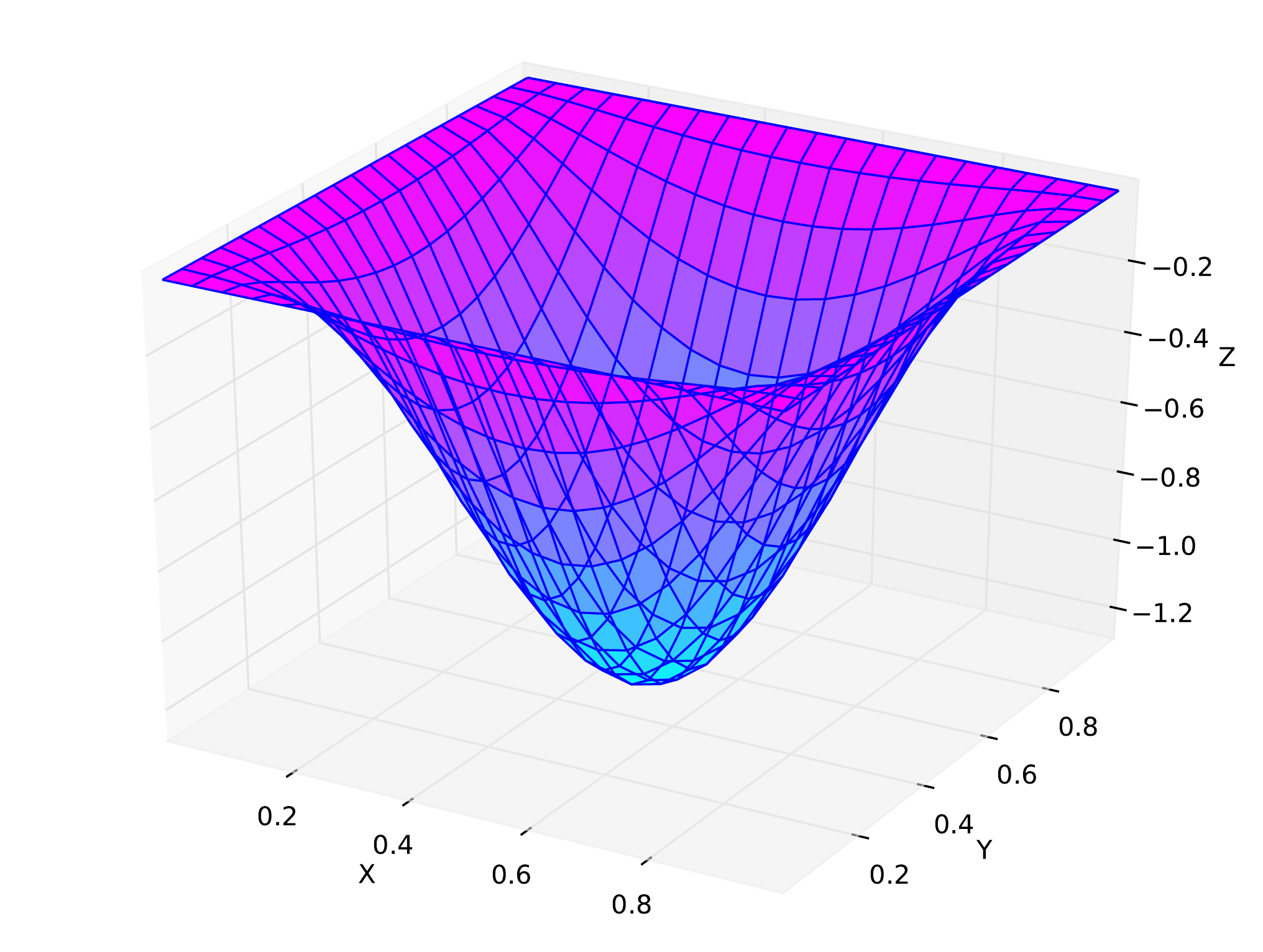}
\caption{Solution of Biharmonic equation (Example 5b), Type II BCs, n=20.}
\label{fig:ex5b}
\end{figure}

Biharmonic equation that is solved in this example has Type II boundary conditions and cannot be reformulated as a system of two coupled Poisson equations. This example is also taken from \cite{Yakhot}.  We will consider two cases, one having the exact solution\\
\textit{Case 5a:}
\[
f(x,y) = 56400(a^2-10ax+15x^2)(b-y)^2y^4
\]\[
+18800x^2(6a^2-20ax+15x^2)y^2(6b^2-20by+15y^2)  
\]
\begin{equation}
+56400(a-x)^2x^4(b^2-10by+15y^2),
\end{equation}
And one where the exact solution is unknown
\textit{Case 5b:}
\begin{equation}
f(x,y) = p_0 = const,
\end{equation}
where $p_0$ takes value of 1000. The notation has physical significance, Example 5b represents the case of a plate clamped at all four sides, and exposed to uniform load of fluid pressure.
Example 5a admits the exact solution
\begin{equation}
u_{exact}(x,y) =2350 x^4 (x-a)^2 y^4 (y-b)^2.
\end{equation}
In the assembling procedure, instead of writing the Eq. \ref{eq:a_assemble} for the nodes laying right next to the egde nodes, we use Eq. \ref{eq:b_assemble}, which to remind once again, uses $(x,y)$ values of the points on the boundary edge. These edge points are all immediate neighbours of those collocations points, we would normally write the equation for. We note that the Neumann boundary conditions are not written for the corner nodes. Numerical results, consistent with previous examples, is shown in Fig. \ref{fig:ex5a}, Fig. \ref{fig:ex5aabserr} and Table \ref{l2err5a}.
When there is no exact solution, we need to set up a criteria what solution to except as converged one. In all previous examples, which allowed exact solution, we see that the high-order accuracy is achieved by a comparatively small number of nodal points. For 21 nodal point in each direction (order of polynomial n=20), the $L^2$ relative error norm is of the order $10^{10}$ to $10^{11}$. We should have the additional confidence in the result, if the spatial variation of solution expressed in terms of local maxima and minima, is not significant within the domain. Fig. \ref{fig:ex5b} shows numerical solution with Bernstein polynomials or order n=20 in each direction. 

\begin{table}
\centering
\begin{tabular} {|c|c|}
\hline
n & $L^2$ rel. error \\ \hline \hline
8 & $2.506 \times 10^{-14}$ \\
10 & $1.064 \times 10^{-14}$ \\
12 & $3.773 \times 10^{-13}$ \\
14 & $1.174 \times 10^{-11}$ \\
20 & $3.853 \times 10^{-10}$ \\
\hline
\end{tabular}
\caption{$L^2$ relative error norm for the Example 5a.}
\label{l2err5a}
\end{table}

\section{Conclusions}
In this paper a novel formulation of the collocation method using Bernstein polynomials is proposed. The main reason why the collocation method is chosen, are it's flexibility and simple implementation.
The methodology presented in this paper has been implemented in \textsf{bernstein-poly} code \cite{b-poly-code}, and several examples have been shown, where the elliptic boundary value problems in two dimensional domains were succesfully solved. Numerical results obtained by this method were compared with existing exact solutions. Excellent agreement and high accuracy is achievent even with small number of basis polynomials. \\
Trough extstensive testing we have concluded that the three components of the algorithm used here: defining polynomials over general interval, the non-recursive formulation for derivatives and the use of multiplicative formula for individual binomial coefficients significantly enhance capabilities of the present procedure related to the previous Bernstein polynomial methods in terms of the flexibility and speed.\\
As we have seen in log-log plots of $L^2$ relative error norm in all examples, the error decreases exponentially as the order of polynomial increases, until, around $n=17-20$, when it changes the character and continually increases for higher values of n. We noticed the same character of error variation with n in all the examples above, which is related to the loss in accuracy in floating-point arithmetics with large numbers originating from factorials in the definition of Bernstein polynomials. This suggest further direction of investigation - developing a method based on a principle of domain decomposition. In that case in each "element" we would keep moderate degree in Bernstein polynomial, and the method would resemble spectral element method. Similar has been done e.q. \cite{Kopriva}, where Chebyshev multidomain method has been put forward. That study has served as a basis for development Spectral Difference Method. 
Another advantage of the domain decomosition approach is that it produces sparse matrices, with the sparsity pattern usual for tensor product grid discretizations. \\
Future plan is to develop \textsf{bernstein-poly} to the point where highly accurate solution of Navier-Stokes equations in complex three-dimensional domains is possible.


\newpage

\appendix

\section*{Appendix A. Implementation of a new BVP in \textsf{bernstein-poly}}

Setting up a new elliptic boundary value problem in \textsf{bernstein-poly} code \cite{b-poly-code} is briefly described here on the case of Poisson equation presented in the Example 3.
Main program is located in \textsf{collocation\_test\_2D.py}, where the user needs to define solution domain by specifying it's x and y axis extents, and the order of approximating polynomials - n.

\begin{lstlisting}[label=code1,caption=Setting up domain borders and approximating polynomial order.]
# Bernstein polynomial order - the same order in both directions.
n = 12 
m = 12
# The number of unknowns - all nodes counted for the case with non-homogenous BCs.
nvar = (n+1)*(m+1) 
...
# Solution interval [x1,x2] x [y1,y2].
x1 = 0.   
x2 = 1. 

y1 = 0.
y2 = 1. 
...
# Uniform mesh.
nd = linspace(x1,x2,n+1) 
\end{lstlisting}

Specifying the problem generaly described by Eq. (\ref{eq:L-operator}) is straight-forward:
\begin{lstlisting}[label=code2,caption=Defining RHS vector entry]
def rhs(x,y):
    return -(6*x*y*(1-y)-2*x**3)
\end{lstlisting}

\begin{lstlisting}[label=code3,caption=Defining system matrix entry]
def lhs(i,j,n,m,x1,x2,y1,y2,x,y):
    return -laplacian(i,j,n,m,x1,x2,y1,y2,x,y) 
\end{lstlisting}

\begin{lstlisting}[label=code4,caption=Defining the exact solution.]
def exact_solution(x,y):
	return y*(1-y)*x**3
\end{lstlisting}

Boundary conditions have to be specified next. First we'll have a look how does the LHS matrix and RHS vector assembly look like (Listings 5 and 6).

\begin{lstlisting}[label=code5,caption=Forming RHS vector in the main function.]
f = zeros(nvar) # RHS vector 
for i in range(0,n+1):
    x = nd[i]
    for j in range(0,m+1):
        y = nd[j]
        node = (m+1)*(i-1)+j  
#       Non-homogenous BCs...
        if (i==0):
#       Left side:: Run trough all betas
	    f[node] = bc_left(x,y)    
        elif (j==0):
#       Bottom:: Run trough all betas
	    f[node] = bc_bottom(x,y) 
        elif (i==n):
#       Right: Run trough all betas
	    f[node] = bc_right(x,y) 
        elif (j==n):
#       Top: Run trough all betas
	    f[node] = bc_top(x,y) 
        else:  
	    f[node] = rhs(x,y) 
\end{lstlisting}

\begin{lstlisting}[label=code6,caption=Forming system matrix]
K = zeros( (nvar,nvar) ) # LHS matrix 
for i in range(0,n+1):
    x = nd[i]
    for j in range(0,m+1):
        y = nd[j]
        node = (m+1)*(i-1)+j # node defines specific location on a grid
#       Non-homogenous BCs...
        if (i==0):
#       Left side:: Run trough all betas
             for k in range(0,n+1):
                for l in range(0,m+1):
                    jfun = (m+1)*(k-1)+l
	            K[node,jfun] = Dirichlet(l,k,n,m,x1,x2,y1,y2,x,y)    
        elif (j==0):
#       Bottom:: Run trough all betas
            for k in range(0,n+1):
                for l in range(0,m+1):
                    jfun = (m+1)*(k-1)+l
	            K[node,jfun] = Dirichlet(l,k,n,m,x1,x2,y1,y2,x,y) 
        elif (i==n):
#       Right: Run trough all betas
            for k in range(0,n+1):
                for l in range(0,m+1):
                    jfun = (m+1)*(k-1)+l
	            K[node,jfun] = Dirichlet(l,k,n,m,x1,x2,y1,y2,x,y) 
        elif (j==n):
#       Top: Run trough all betas
            for k in range(0,n+1):
                for l in range(0,m+1):
                    jfun = (m+1)*(k-1)+l
	            K[node,jfun] = Dirichlet(l,k,n,m,x1,x2,y1,y2,x,y) 
        else:
#       Interior: Run trough all betas
            for k in range(0,n+1):
                for l in range(0,m+1):
                    jfun = (m+1)*(k-1)+l
	            K[node,jfun] = lhs(l,k,n,m,x1,x2,y1,y2,x,y)  
\end{lstlisting}
	
\begin{lstlisting}[label=code7,caption=System matrix entry originating from the boundary operator (\ref{eq:B-operator})]    
def Dirichlet(i,j,n,m,x1,x2,y1,y2,x,y):
# Matrix entry for Dirichlet BCs:
    return basis_fun_eval(i,n,x1,x2,x)*basis_fun_eval(j,m,y1,y2,y)

def Neumann_y_dir(i,j,n,m,x1,x2,y1,y2,x,y):
# Matrix entry for Neumann BCs:
    return basis_fun_eval(i,n,x1,x2,x)*basis_fun_der(1,j,m,y1,y2,y) 
\end{lstlisting}     

RHS vector entry for the collocation point located at domain boundary is evaluated according to function $g(x)$ Eq. (\ref{eq:B-operator}).
\begin{lstlisting}[label=code8,caption=Evaluating RHS vector entry for the points on domain boundary.]
def bc_right(x,y):
# Boundary condition for right edge of the rectangle
    return y*(1-y)
\end{lstlisting}     

The beauty of the proposed method and it's Python implementation, is that non-trivial problems descibed using PDEs are solved in a very small number of command lines. The guiding principle behind the development of the code is modularity and escalation towards solution of problems with greater complexity. Readers are encouraged to use and upgrade \textsf{bernstein-poly} in their own research. 
 

\begin{thebibliography}{1}

\bibitem{Lorentz}
G.~G. Lorentz.
\newblock {\em Bernstein Polynomials}.
\newblock University of Toronto Press, Toronto, Canada, 1953.

\bibitem{Bhatti2007}
M.~I. Bhatti and P.~Bracken.
\newblock Solutions of differential equations in a {B}ernstein polynomial
  basis.
\newblock {\em Journal of Computational and Applied mathematics},
  205(1):272--280, 2007.
  
\bibitem{Yousefi}
S.A.~Yousefi, Z.~Barikbin, M.~Denhgan.
\newblock Bernstein Ritz-Galerkin method for solving an initial-boundary value
  problem that combines neumann and integral condition for the wave equation.
\newblock {\em Numerical Methods for Partial Differential Equations},
  26:1236--1246, 2009.
  
\bibitem{Doha}
E. H.~Doha, A. H.~Bhrawy, M.~A.~Saker.
\newblock On the derivatives of {B}ernstein polynomials: An application for
  the solution of high-even-order deifferential equations.
\newblock {\em Boundary Value Problems}, 2011.

\bibitem{Boyd}
J.~Boyd.
\newblock {\em {C}hebyshev and {F}ourier {S}pectral {M}ethods}.
\newblock Dover, 2nd edition, 2001.

\bibitem{Fornberg}
B.~Fornberg.
\newblock {\em A Practical Guide to Pseudospectral Methods}.
\newblock Cambridge University Press, 1998.

\bibitem{ICASEref}
M.~Y.Hussaini, C.~L.Streett, T.~A.Zang.
\newblock Spectral Methods for Partial Differential Equations.
\newblock NASA Contractor Report, NASA-CR-172248, 1983.

\bibitem{transformation}
R.~T. Farouki.
\newblock Legendre-Bernstein basis transformations.
\newblock {\em J. Comput. Appl. Math.}
\newblock 119: 145--160, 2000.

\bibitem{b-poly-code}
\newblock {http://code.google.com/p/bernstein-poly/}

\bibitem{Hon}
Y.~C.Hon, W.~Chen.
\newblock Boundary knot method for 2D and 3D {H}elmholtz and convection-diffusion problems under complicated geometry.
\newblock {\em International Journal of Numerical Methods in Engineering},
56:1931--1948, 2003.

\bibitem{Yakhot}
M.~Arad, A.~Yakhot, G.~Ben-Dor. 
\newblock A Highly Accurate Numerical Solution of a Biharmonic Equation.
\newblock {\em Numerical Methods for Partial Differential Equations},
  13:375--391, 1997.

\bibitem{Kopriva}
D.~A.Kopriva, J.~H.Kolias.
\newblock A conservative staggered-grid Chebyshev multidomain method for compressible flows.
\newblock {\em Journal of Computational Physics},
125(1):244--261, 1996.


\end{thebibliography}






\end{document}